\documentclass[reqno,centertags, draft]{amsart}
\usepackage{amsmath,amsthm,amscd,amssymb,latexsym,upref,enumerate}

\newcommand{\bbC}{{\mathbb C}}
\newcommand{\bbD}{{\mathbb D}}

\newcommand{\bbN}{{\mathbb N}}

\newcommand{\bbQ}{{\mathbb Q}}
\newcommand{\bbR}{{\mathbb R}}

\newcommand{\bbZ}{{\mathbb Z}}

\newcommand{\cE}{{\mathcal E}}

\newcommand{\cM}{{\mathcal M}}
\newcommand{\cN}{{\mathcal N}}

\newcommand{\cR}{{\mathcal R}}

\newcommand{\loc}{\text{\rm{loc}}}

\newcommand{\rank}{\text{\rm{rank}}}

\newcommand{\Arg}{\text{\rm{Arg}}}
\newcommand{\Arc}{\text{\rm{Arc}}}
\newcommand{\dom}{\text{\rm{dom}}}
\newcommand{\tr}{\text{\rm{tr}}}

\newcommand{\ess}{\text{\rm{ess}}}
\newcommand{\ac}{\text{\rm{ac}}}
\renewcommand{\sc}{\text{\rm{sc}}}
\newcommand{\s}{\text{\rm{s}}}
\newcommand{\pp}{\text{\rm{pp}}}

\newcommand{\supp}{\text{\rm{supp}}}

\newcommand{\id}{\text{\rm{id}}}

\newcommand{\Om}{\Omega}
\newcommand{\om}{\omega}
\newcommand{\si}{\sigma}
\newcommand{\la}{\lambda}
\newcommand{\La}{\Lambda}
\newcommand{\al}{\alpha}

\newcommand{\de}{\delta}
\newcommand{\te}{\theta}
\newcommand{\ze}{\zeta}

\newcommand{\eps}{\varepsilon}

\newcommand{\bi}{\bibitem}
\newcommand{\no}{\notag}
\newcommand{\lb}{\label}
\newcommand{\f}{\frac}

\newcommand{\ol}{\overline}

\newcommand{\bs}{\backslash}

\newcommand{\abs}[1]{\left\lvert#1\right\rvert}

\newcommand{\dott}{\,\cdot\,}
\newcommand{\dD}{{\partial\mathbb{D}}}

\newcommand{\Cl}{\mathbb{C}_{\ell}}
\newcommand{\Cr}{\mathbb{C}_{r}}

\newcommand{\st}{\,|\,}
\newcommand{\ltz}{{\ell^2(\bbZ)}}
\newcommand{\lt}[1]{{\ell^2(#1)}}

\newcommand{\deven}{\delta_{\rm even}}
\newcommand{\dodd}{\delta_{\rm odd}}

\renewcommand{\Re}{\text{\rm Re}}
\renewcommand{\Im}{\text{\rm Im}}
\renewcommand{\ln}{\text{\rm ln}}

\allowdisplaybreaks
\numberwithin{equation}{section}

\newtheorem{theorem}{Theorem}[section]

\newtheorem{lemma}[theorem]{Lemma}

\newtheorem{definition}[theorem]{Definition}
\newtheorem{example}[theorem]{Example}
\theoremstyle{definition}
\newtheorem{remark}[theorem]{Remark}

\begin{document}

\title[Essential closures and ac spectra for reflectionless operators]
{Essential Closures and AC Spectra for Reflectionless CMV, Jacobi, and 
Schr\"odinger Operators Revisited}
\author[F.\ Gesztesy, K.\ A.\ Makarov, and M.\ Zinchenko]
{Fritz Gesztesy, Konstantin A.\ Makarov, and Maxim Zinchenko}

\address{Department of Mathematics,
University of Missouri, Columbia, MO 65211, USA}
\email{fritz@math.missouri.edu}
\urladdr{http://www.math.missouri.edu/personnel/faculty/gesztesyf.html}

\address{Department of Mathematics, University of
Missouri, Columbia, MO 65211, USA}
\email{makarov@math.missouri.edu}
\urladdr{http://www.math.missouri.edu/personnel/faculty/makarovk.html}

\address{Department of Mathematics,
California Institute of Technology, Pasadena, CA 91125, USA}
\email{maxim@caltech.edu}
\urladdr{http://www.math.caltech.edu/\~{}maxim}

\date{\today}
\dedicatory{Dedicated with great pleasure to L.\ J.\ Lange on the occasion of his 80th birthday.}

\begin{abstract}
We provide a concise, yet fairly complete discussion of the concept of essential closures of subsets of the real axis and their intimate connection with the topological support of absolutely continuous measures.

As an elementary application of the notion of the essential closure of subsets of $\bbR$
we revisit the fact that CMV, Jacobi, and Schr\"odinger operators, reflectionless on a set $\cE$ of positive Lebesgue measure, have absolutely continuous spectrum on the essential closure ${\ol \cE}^e$ of the set $\cE$ (with uniform multiplicity two on $\cE$). Though this result in the case of Schr\"odinger and Jacobi operators is known to experts, we feel it nicely illustrates the concept and usefulness of essential closures in the spectral theory of classes of reflectionless differential and difference operators.   
\end{abstract}

\maketitle

\section{Introduction}\lb{s1}

In this note we revisit the notion of essential closures of subsets of the real line and their intimate connection with the topological support of absolutely continuous measures.  
As an elementary application of this concept we consider unitary CMV operators and self-adjoint Jacobi operators on $\bbZ$, and Schr\"odinger operators on $\bbR$, which are reflectionless on a set $\cE$ of positive Lebesgue measure and recall the elementary proof that the essential closure of $\cE$, denoted by $\ol{\cE}^e$, belongs to their  absolutely continuous spectrum. 

We emphasize that this paper is in part of expository nature. Still, we feel it is worth the effort to systematically highlight properties of essential closures of sets, and recall how this concept naturally leads to the existence of absolutely continuous spectra of certain well-known classes of operators, especially, in the presence of a reflectionless condition.  

While our emphasis here is on reproving the existence of absolutely continuous spectrum with most elementary methods, the absence of singular spectrum is quite a distinct matter that typically requires entirely different methods not discussed in this paper (in the context of reflectionless operators, see however, \cite{GY06}, \cite{GZ08}, 
\cite{Ko84}, \cite{KK88}, \cite{PY06}, \cite{SY95a}--\cite{SY97}, and the literature cited therein). 

Next, we briefly single out Schr\"odinger operators and illustrate the notion of being reflectionless: {\it Reflectionless} (self-adjoint) 
Schr\"odinger operators $H$ in $L^2(\bbR;dx)$ can be characterized, for instance, by the fact that for all $x\in\bbR$ and for a.e.\ $\lambda\in \sigma_{\rm ess}(H)$, the diagonal Green's function of $H$ has purely imaginary normal boundary values,
\begin{equation}
G(\lambda+i0,x,x) \in i\bbR.    \lb{1.3}
\end{equation}
Here $\sigma_{\rm ess}(H)$ denotes the essential spectrum of $H$ (we assume
$\sigma_{\rm ess}(H)\neq\emptyset$) and
\begin{equation}
G(z,x,x')=(H-zI)^{-1}(x,x'), \quad z\in\bbC\backslash\sigma(H), \lb{1.4}
\end{equation}
denotes the integral kernel of the resolvent of $H$. This global notion of reflectionless Schr\"odinger operators can of course be localized and extends to subsets of 
$\sigma_{\rm ess}(H)$ of positive Lebesgue measure. In the actual body of our paper we will use an alternative definition of the notion of reflectionless Schr\"odinger operators conveniently formulated directly in terms of half-line Weyl--Titchmarsh functions, we refer to Definitions \ref{d3.2}, \ref{d2.2}, and \ref{d4.2} for more details. For various discussions of classes of reflectionless differential and difference operators, we refer, 
for instance, to Craig \cite{Cr89}, De Concini and Johnson \cite{DJ87}, 
Deift and Simon \cite{DS83}, Gesztesy, Krishna, and Teschl \cite{GKT96}, Gesztesy and Yuditskii \cite{GY06}, Johnson \cite{Jo82}, Kotani \cite{Ko84}--\cite{Ko87b}, Kotani and Krishna \cite{KK88}, Peherstorfer and Yuditskii \cite{PY06}, \cite{PY07}, Remling \cite{Re07}, \cite{Re08}, Sims \cite{Si07a}, Sodin and Yuditskii \cite{SY95a}--\cite{SY96}, 
and Vinnikov and Yuditskii \cite{VY02}. In particular, we draw attention to two  recent papers by Remling \cite{Re07}, \cite{Re08}, that illustrate in depth the ramifications of the existence of absolutely continuous spectra in one-dimensional problems. 

Analogous considerations apply to Jacobi operators (see, e.g., \cite{CL90}, 
\cite{Te00} and the literature cited therein) and CMV operators (see 
\cite{Si04a}--\cite{Si05a}, \cite{Si07} and the extensive list of references provided therein and \cite{GZ06a} for the notion of reflectionless CMV operators). For an exhaustive list of references on reflectionless Jacobi and Schr\"odinger operators we refer to the bibliography in \cite{GZ08}. 

In Section \ref{s2} we review basic facts on essential closures of sets and essential supports of measures. In Section \ref{s3} we consider CMV, Jacobi, and 
Schr\"odinger operators reflectionless on sets $\cE$ of positive Lebesgue measure and recall that the essential closure $\ol{\cE}^e$ of $\cE$ belongs to their absolutely continuous spectrum. For brevity, we provide proofs in the CMV case only as this case has received considerably less attention when compared to Jacobi and 
Schr\"odinger operators.

The methods employed in Section \ref{s3} are elementary and based on the facts discussed in Section \ref{s2} and on the material presented in Appendices 
\ref{sB} and \ref{sC}. The latter provide a nutshell-type treatment of properties of Herglotz and Caratheodory functions as well as certain elements of Weyl--Titchmarsh and spectral multiplicity theory for self-adjoint Jacobi and Schr\"odinger operators on 
$\bbZ$ and $\bbR$, and unitary CMV operators on $\bbZ$.

\section{Basic facts on essential closures of sets \\
and essential supports of measures} \lb{s2}

The following material on essential closures of subsets of the real line and the unit circle and  essential supports of measures is well-known to experts, but since no
comprehensive treatment in the literature appears to exist in one place, we have collected the relevant facts in this section.

For basic facts on measures on $\bbR$ relevant to this section we
refer, for instance, to \cite{Ar57}--\cite{AD64}, \cite[p.\
179]{CFKS87}, \cite{DSS94}, \cite{Gi84}--\cite{GP87}, \cite[Sect.\
V.12]{PF92}, \cite[p.\ 140--141]{RS78}, \cite{Si95}, \cite{SW86}.
All measures in this section will be assumed to be nonnegative
without explicitly stressing this fact again.

Since Borel and Borel--Stieltjes measures are incomplete (i.e., not
any subset of a set of measure zero is measurable) we will enlarge
the Borel $\sigma$-algebra to obtain the complete Lebesgue and
Lebesgue--Stieltjes measures. We recall the standard Lebesgue
decomposition of a measure $d\mu$ on $\bbR$ with respect to Lebesgue
measure $dx$ on $\bbR$,
\begin{align}
&d\mu=d\mu_{\ac}+d\mu_{\rm s}=d\mu_{\ac}+d\mu_{\sc} +d\mu_{\rm
pp}, \lb{A.1} \\
&d\mu_{\ac}=fdx, \quad 0\leq f\in L^1_{\rm loc}(\bbR;dx), \lb{A.2}
\end{align}
where $d\mu_{\ac}$, $d\mu_{\s}$, $d\mu_{\sc}$, and $d\mu_{\pp}$
denote the absolutely continuous, singular, singularly continuous,
and pure point parts of $d\mu$, respectively.

In the following, the Lebesgue measure of a Lebesgue measurable set
$S\subseteq\bbR$ will be denoted by $|S|$ and all sets whose
$\mu$-measure or Lebesgue measure is considered are always assumed
to be Lebesgue--Stieltjes or Lebesgue measurable, etc.

\begin{definition} \lb{dA.1}
Let $d\mu$ be a Lebesgue--Stieltjes measure and suppose $S$ and $S'$
are $\mu$-measurable. \\
$(i)$ $S$ is called a {\it support} of $d\mu$ if
$\mu(\bbR\bs S)=0$. \\
$(ii)$ The smallest closed support of $d\mu$ is called the {\it
topological support} of $d\mu$ and denoted by $\supp \, (d\mu)$. \\
$(iii)$ $S$ is called an {\it essential} $($or {\it minimal\,$)$
support} of $d\mu$ $($relative to Lebesgue measure $dx$ on $\bbR$$)$ if
$\mu(\bbR\bs S)=0$, and $S'\subseteq S$ with $S'$
$|\cdot|$-measurable, $\mu(S')=0$ imply $|S'|=0$.
\end{definition}

\begin{remark} \lb{rA.2}
Item $(iii)$ in Definition \ref{dA.1} is equivalent to \\
$(iii')$ $S$ is called an {\it essential} $($or {\it minimal\,$)$
support} of $d\mu$ (relative to Lebesgue measure $dx$ on $\bbR$) if
$\mu(\bbR\bs S)=0$, and $S'\subseteq S$, $\mu(\bbR\bs S')=0$ imply
$|S\bs S'|=0$.
\end{remark}

\begin{lemma} [\cite{Gi84}]\lb{lA.3}
Let $S, S' \subseteq\bbR$ be $\mu$- and $|\cdot|$-measurable. Define
the relation $\sim$ by $S\sim S'$ if
\begin{equation}
\mu(S\Delta S')=|S\Delta S'|=0  \lb{A.3}
\end{equation}
$($where $S\Delta S'=(S\bs S')\cup (S'\bs S)$$)$. Then $\sim$ is an
equivalence relation. Moreover, the set of all essential supports of
$d\mu$ is an equivalence class under $\sim$.
\end{lemma}

\begin{example} \lb{eA.4}
Let $d\mu_{\pp}$ be a finite pure point measure and
\begin{equation}
\mu_{\pp}(\{x\})=\begin{cases} c(x)>0, & x\in [0,1]\cap \bbQ, \\
0, & \text{otherwise}.
\end{cases}
\lb{A.4}
\end{equation}
Then,
\begin{equation}
\supp \, (d\mu_{\pp})=[0,1]. \lb{A.5}
\end{equation}
However, since $[0,1]\cap \bbQ$ is an essential support of
$d\mu_{\pp}$ and since $|[0,1]\cap \bbQ|=0$, also
\begin{equation}
|S_{\mu_{\pp}}|=0  \lb{A.6}
\end{equation}
for any other essential support $S_{\mu_{\pp}}$ of $d\mu_{\pp}$.
\end{example}

\begin{remark} \lb{rA.5} ${}$ \\ 
$(i)$ Any two essential supports of $d\mu$ differ at most
by sets of Lebesgue measure zero. \\
$(ii)$ Assume $d\mu =d\mu_{\ac}$ and let $S, S' \subseteq\bbR$ 
be $\mu$- and $|\cdot|$-measurable. Then $|S\Delta S'|=0$ implies 
$\mu(S\Delta S')=0$. In particular, $|S\Delta S'|=0$ implies that $S$ is 
an essential support of $d\mu$ if and only if $S'$ is. \\
Indeed, one can use the following elementary relations, 
\begin{align}
&S_1=(S_1\cap S_2)\cup (S_1\bs S_2), \quad
S_2=(S_2\cap S_1)\cup (S_2\bs S_1),  \lb{A.7} \\
&S_1\cup(S_2\bs S_1)=S_2\cup(S_1\bs S_2),  \lb{A.8}
\end{align}
valid for any subsets $S_j\subseteq\bbR$, $j=1,2$.
\end{remark}

\begin{definition} \lb{dA.6}
Let $A\subseteq\bbR$ be Lebesgue measurable. Then the {\it essential
closure} ${\ol A}^e$ of $A$ is defined as
\begin{equation}
{\ol A}^e =\{x\in\bbR\,|\, \text{for all $\eps>0$$:$} \, 
|(x-\eps,x+\eps)\cap A|>0\}. \lb{A.9}
\end{equation}
\end{definition}

The following is an immediate consequence of Definition \ref{dA.6}.

\begin{lemma} \lb{lA.7}
Let $A, B, C \subseteq\bbR$ be Lebesgue measurable. Then,
\begin{align}
& \text{$(i)$ \; If $A\subseteq B$ then ${\ol A}^e\subseteq {\ol B}^e$.}
\lb{A.10A} \\
& \text{$(ii)$\;\;If $|A|=0$ then ${\ol A}^e=\emptyset$.}  \lb{A.10} \\
& \text{$(iii)$\;If $A=B\cup C$ with $|C|=0$, then ${\ol A}^e={\ol B}^e$.}
\lb{A.10B} \\
& (iv)\,\,\, \ol{A}^e \text{ is a closed set.}  \lb{A.10a} 
\end{align}
\end{lemma}
\begin{proof} 
Since items $(i)$--$(iii)$ are obvious, it suffices to focus on item $(iv)$. We will 
show that the set
\begin{equation}
\bbR\bs{\ol A}^e=\{x\in\bbR\,|\, \text{there is an $\eps_0>0$ such
that } |(x-\eps_0,x+\eps_0)\cap A|=0\} \lb{A.11}
\end{equation}
is open. Pick $x_0\in\bbR\bs{\ol A}^e$, then there is an $\eps_0>0$
such that $|(x_0-\eps_0,x_0+\eps_0)\cap A|=0$. Consider
$x_1\in(x_0-(\eps_0/2),x_0+(\eps_0/2))$ and the open ball
$S(x_1;\eps_0/2)$ centered at $x_1$ with radius $\eps_0/2$. Then,
\begin{equation}
|S(x_1;\eps_0/2)\cap A|\leq |(x_0-\eps_0,x_0+\eps_0)\cap A|=0
\lb{A.12}
\end{equation}
and hence $x_1\in\bbR\bs{\ol A}^e$ and $S(x_0;\eps_0/2)\subseteq
\bbR\bs {\ol A}^e$. Thus,
$\bbR\bs{\ol A}^e$ is open. 
\end{proof}

\begin{example} \lb{eA.8} ${}$ \\
$(i)$ Consider $d\mu_{\pp}$ in Example \ref{eA.4}. Let $S_{\mu_{\rm
pp}}$ be any essential support of $d\mu_{\pp}$. Then
${\ol {S_{\mu_{\pp}}}}^e=\emptyset$ by \eqref{A.10}. \\
$(ii)$ Consider $A=[0,1]\cup\{2\}$. Then ${\ol A}^e =[0,1]$.
\end{example}

\begin{lemma} \lb{lA.10}
Let $S_1$ and $S_2$ be essential supports of $d\mu$. Then,
\begin{equation}
{\ol{S_1}}^e = {\ol{S_2}}^e.  \lb{A.13}
\end{equation}
\end{lemma}
\begin{proof}
Since $|S_1\bs S_2|=|S_2\bs S_1|=0$, \eqref{A.13} follows from
\eqref{A.7} and \eqref{A.10B}.
\end{proof}

Actually, one also has the following result: 

\begin{lemma} \lb{lA.11}
Let $d\mu=d\mu_{\ac}=fdx$, $0\leq f\in L^1_{\rm loc}(\bbR)$. If $S$
is any essential support of $d\mu$, then,
\begin{equation}
{\ol S}^e = \ol{\{x\in\bbR\,|\, f(x)>0\}}^e=\supp \, (d\mu).  \lb{A.14}
\end{equation}
\end{lemma}
\begin{proof}
Since $\{x\in\bbR\,|\, f(x)>0\}$ is an essential support of $d\mu$, it
suffices to prove
\begin{equation}
\ol{\{x\in\bbR\,|\, f(x)>0\}}^e=\supp \, (d\mu).  \lb{A.15}
\end{equation}
We denote $U=\bbR\bs \supp \, (d\mu)$. Then $U$ is the largest open
set that satisfies $\mu(U)=0$. Next, let
$U'=\bbR\bs\ol{\{x\in\bbR\,|\, f(x)>0\}}^e$. By Lemma \ref{lA.7}\,$(iv)$, 
$U'$ is open. \\
``$\supseteq$'': Let $x\in U'$. Then there is an $\eps_0>0$ such
that
\begin{equation}
|(x-\eps_0,x+\eps_0)\cap \{y\in\bbR\,|\, f(y)>0\}|=0. \lb{A.16}
\end{equation}
Hence,
\begin{equation}
f=0 \; |\cdot|\text{-a.e.\ on $(x-\eps_0,x+\eps_0)$} \lb{A.16a}
\end{equation}
and thus, $\mu((x-\eps_0,x+\eps_0))=0$. The collection of all such
open intervals forms an open cover of $U'$. Then selecting a
countable subcover, one arrives at $\mu(U')=0$. Since $U$ is the
largest open set satisfying $\mu(U)=0$, one infers $U'\subseteq U$
and hence
\begin{equation}
\ol {\{x\in\bbR\,|\, f(x)>0\}}^e \supseteq \supp \, (d\mu). \lb{A.17}
\end{equation}
``$\subseteq$'': Fix an $x\in U$. Since $U$ is open, there is an
$\eps_0>0$ such that $(x-\eps_0,x+\eps_0)\cap \supp \,
(d\mu)=\emptyset$. Thus, $\mu((x-\eps_0,x+\eps_0))=0$. Actually,
$\mu(B)=0$ for all $\mu$-measurable $B\subseteq (x-\eps_0,x+\eps_0)$
and hence $f=0$ $|\cdot|$-a.e.\ on $(x-\eps_0,x+\eps_0)$. Consequently, one
obtains
\begin{equation}
|(x-\eps_0,x+\eps_0)\cap \{y\in\bbR\,|\, f(y)>0\}|=0 \lb{A.18}
\end{equation}
which implies $x\in U'$, $U'\supseteq U$, and
\begin{equation}
\ol {\{x\in\bbR\,|\, f(x)>0\}}^e \subseteq \supp \, (d\mu). \lb{A.19}
\end{equation}
\end{proof}

We remark that a result of the type \eqref{A.14} has been noted
in \cite[Corollary 11.11]{Bu97} in the context of general ordinary
differential operators and their associated Weyl--Titchmarsh matrices.
In this connection we also refer to \cite[p.\ 301]{Te00} for a
corresponding result in connection with Herglotz functions and their
associated measures.

\begin{lemma} \lb{lA.9}
Let $A\subseteq\bbR$ be Lebesgue measurable. Then,
\begin{align} 
& (i)\quad \,\, \ol{A}^e \subseteq \ol{A}.  \lb{A.10b}  \\
&(ii)\quad \big|A\bs\ol{A}^e\big|=0. \lb{A.26} \\
& (iii)\quad \!\!\ol{\big({\ol A}^e\big)}^e = {\ol A}^e.
\lb{A.10c} \\ 
&(iv)\!\quad \big|\ol{A}^e\big|\geq|A|. \lb{A.28}
\end{align}
\end{lemma}
\begin{proof}
$(i)$ Let $x\in {\ol A}^e$. Then for all $\eps>0$,
$|(x-\eps,x+\eps)\cap A|>0$. Choose $\eps_n=1/n$, $n\in\bbN$, then
$(x-\eps_n,x+\eps_n)\cap A \neq \emptyset$ and we may choose an
$x_n\in (x-\eps_n,x+\eps_n)\cap A$. Since $x_n\to x$ as $n\to
\infty$, $x\in \ol A$ and hence
${\ol A}^e\subseteq \ol A$. \\
$(ii)$ Let $f=\chi_A$ be the characteristic function of the set $A$. Then
Lemma \ref{lA.11} applied to the measure $d\mu=fdx$ implies that
$\ol{A}^e = \supp\,(d\mu)$. Thus, \eqref{A.26} follows from
Definition \ref{dA.1}\,$(i)$,
\begin{equation}
0 = \mu(A\bs\supp\,(d\mu)) = |A\bs\supp\,(d\mu)| = \big|A\bs\ol{A}^e\big|.
\end{equation}
$(iii)$ By $(i)$, ${\ol A}^e\subseteq \ol A$. Hence,
$\ol{\big({\ol A}^e\big)}^e \subseteq \ol{{\ol A}^e}={\ol A}^e$ since
${\ol A}^e$ is closed by Lemma \ref{lA.7}\,$(iv)$. Conversely, 
let $A_1 = A\cap\ol{A}^e$. Then \eqref{A.26} and Lemma \ref{lA.7}\,$(iii)$ 
yield
\begin{equation}
\ol{A}^e=\ol{A_1\cup\big(A\bs\ol{A}^e\big)}^e = \ol{A_1}^e. \label{A.31}
\end{equation}
Since $A_1\subseteq\ol{A}^e=\ol{A_1}^e$, it follows from 
Lemma \ref{lA.7}\,$(i)$ that $\ol{A_1}^e\subseteq\ol{\big(\ol{A_1}^e\big)}^e$.
Thus, \eqref{A.31} implies $\ol{A}^e\subseteq\ol{\big(\ol{A}^e\big)}^e$, and
hence \eqref{A.10c} holds. \\ 
$(iv)$ Equation \eqref{A.26} and 
$A=\big(A\cap\ol{A}^e\big)\cup\big(A\bs\ol{A}^e\big)$ imply \eqref{A.28},
\begin{equation}
|A| = \big|A\cap\ol{A}^e\big| \leq \big|\ol{A}^e\big|.
\end{equation}
\end{proof}

As the following example shows, the inequality in \eqref{A.28} can be strict. 

\begin{example} \lb{eA.12}
Let $\{r_n\}_{n\in\bbN}$ be an enumeration of the rational numbers in $[0,1]$. 
Then the set
\begin{equation}
A=\bigcup_{n\in\bbN} \bigg(r_n-\f{1}{4^n},r_n+\f{1}{4^n}\bigg)
\end{equation}
satisfies 
\begin{equation}
|A|\leq \sum_{n\in\bbN} \bigg|\bigg(r_n-\f{1}{4^n},r_n+\f{1}{4^n}\bigg)\bigg| 
= \sum_{n\in\bbN} \f{2}{4^n} = \f{2}{3}. 
\end{equation}
Next, taking $x\in [0,1]$ and $\varepsilon >0$, then 
$(x-\varepsilon,x+\varepsilon)\cap [0,1]$ contains at least one rational number $r_m$ for some $m\in\bbN$. Hence,
\begin{equation}
|A\cap (x-\varepsilon,x+\varepsilon)| \geq 
\bigg|\bigg(r_m-\f{1}{4^m}, r_m+\f{1}{4^m}\bigg)\cap (x-\varepsilon,x+\varepsilon)\bigg| 
> 0.
\end{equation}
Thus, $\big|\ol{A}^e\big| \geq |[0,1]| =1 > \f{2}{3} \geq |A|$. 

In addition,  
using $|\ol{A}^e| = |\ol{A}^e \cap A| + |\ol{A}^e \bs A|$, one concludes that 
\begin{equation}
|\ol{A}^e \bs A| = |\ol{A}^e| - |\ol{A}^e \cap A| \geq
|\ol{A}^e| - |A| \geq 1-2/3,
\end{equation} 
and hence $\big|\ol{A}^e\bs A\big| \geq 1/3$, but $\big|A\bs\ol{A}^e\big|=0$ by 
\eqref{A.26}.  
\end{example}

\begin{remark} \lb{rA.13}
Similar definitions and results also hold for sets and measures on
the unit circle, denoted by $\dD$ in the following. In particular, we single out the following
ones for later use: Let $A\subseteq\dD$ be Lebesgue measurable, then
the essential closure of $A$ is defined by
\begin{equation}
{\ol A}^e =\big\{e^{i\te}\in\dD \,\big|\, \text{for all $\eps>0$:}\, 
\big|\Arc \big(e^{i(\te-\eps)},e^{i(\te+\eps)}\big)\cap A\big|>0\big\}, \lb{A.33}
\end{equation}
where we used the notation
\begin{equation}
\Arc \big(e^{i\te_1},e^{i\te_2}\big) = \big\{e^{i\te} \, \big|\, 0\leq \te_1<\te<\te_2\big\},
\quad 0\leq \te_1<\te_2< 2\pi. \lb{A.34}
\end{equation}
Let $d\mu$ be a Lebesgue--Stieltjes measure on $\dD$ and suppose
$S\subseteq\dD$ is $\mu$-measurable, then $S$ is called a {\it
support} of $d\mu$ if $\mu(\dD\bs S)=0$. The smallest closed support
of $d\mu$ is called the {\it topological support} of $d\mu$ and
denoted by $\supp\,(d\mu)$. $S$ is called an {\it essential support}
of $d\mu$ (relative to Lebesgue measure on $\dD$) if $\mu(\dD\bs
S)=0$ and $S'\subseteq S$ with $\mu(S')=0$ imply $|S'|=0$. An
essential support $S$ of an absolutely continuous measure
$d\mu=fd\te$ satisfies the identity
\begin{align}
{\ol S}^e = \ol{\big\{e^{i\te}\in\dD \,\big|\, f(e^{i\te})>0\big\}}^e=\supp\,(d\mu). \lb{A.35}
\end{align}
\end{remark}

\section{AC spectra for CMV, Jacobi, and Schr\"odinger operators reflectionless 
on sets of positive Lebesgue measure} \lb{s3}

In this section we apply the results collected on essential closures of subsets of the unit circle and the real line and essential supports of measures in Section \ref{s2} to determine absolutely continuous spectra of CMV, Jacobi, and Schr\"odinger operators reflectionless on sets of positive Lebesgue measure.

We start with the case of unitary CMV operators reflectionless on subsets of the unit circle of positive Lebesgue measure and treat this case in some detail.

Let $\{\al_n\}_{n\in\bbZ}$ be a complex-valued sequence of
Verblunsky coefficients satisfying
\begin{equation}
\al_n\in\bbD = \{z\in\bbC \st |z|<1\}, \quad n\in\bbZ,
\end{equation}
and denote by $\{\rho_n\}_{n\in\bbZ}$ an auxiliary real-valued
sequence defined by
\begin{align}
\rho_n = \big[1-\abs{\al_n}^2\big]^{1/2}, \quad n\in\bbZ.
\end{align}
Then we introduce the associated unitary CMV operator $U$ in $\ltz$
by its matrix representation in the standard basis of $\ltz$,
\begin{align}
U &= \begin{pmatrix} \ddots &&\hspace*{-8mm}\ddots
&\hspace*{-10mm}\ddots &\hspace*{-12mm}\ddots &\hspace*{-14mm}\ddots
&&& \raisebox{-3mm}[0mm][0mm]{\hspace*{-6mm}{\Huge $0$}}
\\
&0& -\al_{0}\rho_{-1} & -\ol{\al_{-1}}\al_{0} & -\al_{1}\rho_{0} &
\rho_{0}\rho_{1}
\\
&& \rho_{-1}\rho_{0} &\ol{\al_{-1}}\rho_{0} & -\ol{\al_{0}}\al_{1} &
\ol{\al_{0}}\rho_{1} & 0
\\
&&&0& -\al_{2}\rho_{1} & -\ol{\al_{1}}\al_{2} & -\al_{3}\rho_{2} &
\rho_{2}\rho_{3}
\\
&&\raisebox{-4mm}[0mm][0mm]{\hspace*{-6mm}{\Huge $0$}} &&
\rho_{1}\rho_{2} & \ol{\al_{1}}\rho_{2} & -\ol{\al_{2}}\al_{3} &
\ol{\al_{2}}\rho_{3}&0
\\
&&&&&\hspace*{-14mm}\ddots &\hspace*{-14mm}\ddots
&\hspace*{-14mm}\ddots &\hspace*{-8mm}\ddots &\ddots
\end{pmatrix}  \no  \\
&= \rho^- \rho \, \deven \, S^{--} + (\ol{\alpha^-}\rho \, \deven - \alpha^+\rho \, \dodd) S^-
- \ol\alpha\alpha^+   \no \\
& \quad + (\ol\alpha \rho^+ \, \deven - \alpha^{++} \rho^+ \, \dodd) S^+
+ \rho^+ \rho^{++} \, \dodd \, S^{++},      \lb{3.4}
\end{align}
where we use the notation for $ f=\{f(n)\}_{n\in\bbZ}\in\ell^{\infty}(\bbZ)$,
\begin{align}
&(S^\pm f)(n)=f(n\pm 1)= f^\pm(n), \quad n\in\bbZ, \no \\
&S^{++}=(S^+)^+, \; S^{--}=(S^-)^-, \, \text{ etc.}   \lb{2.1}
\end{align}
Here terms of the form $-\ol{\alpha_n} \alpha_{n+1}$
represent the diagonal $(n,n)$-entries, $n\in\bbZ$, in the infinite matrix
$U$, and $\deven$ and $\dodd$ denote the characteristic functions of the even and odd integers,
\begin{equation}
\deven = \chi_{_{2\bbZ}}, \quad \dodd = 1 - \deven = \chi_{_{2\bbZ +1}}.
\end{equation}

Moreover, let $M_{1,1}(z,n)$ denote the diagonal element of the
Cayley transform of $U$, that is,
\begin{align}
M_{1,1}(z,n)=((U+zI)(U-zI)^{-1})(n,n) = \oint_\dD
d\Om_{1,1}(\ze,n)\,\f{\ze+z}{\ze-z},& \no
\\
z\in\bbC\bs\sigma(U), \; n\in\bbZ,&
\end{align}
where $d\Om_{1,1}(\cdot,n)$, $n\in\bbZ$, are scalar-valued
probability measures on $\dD$ (cf.\ \cite[Section 3]{GZ06} for more
details). Since for each $n\in\bbZ$, $M_{1,1}(\cdot,n)$ is a
Caratheodory function (i.e., it maps the open unit disk analytically
to the complex right half-plane),
\begin{equation}
\Xi_{1,1}(\ze,n)=\f{1}{\pi}\lim_{r\uparrow 1}
\Im[\ln(M_{1,1}(r\ze,n))] \, \text{ for a.e.\ $\ze\in\dD$} \lb{3.5}
\end{equation}
is well-defined for each $n\in\bbZ$. In particular, for all
$n\in\bbZ$,
\begin{equation}
-1/2 \leq \Xi_{1,1}(\ze,n) \leq 1/2 \, \text{ for a.e.\ $\ze\in\dD$}
\end{equation}
(cf. \cite[Section 2]{GZ06a} for more details).

In the following we will frequently use the convenient abbreviation
$h(\ze)=\lim_{r\uparrow 1} h(r\ze)$, $\ze\in\dD$, whenever the limit is well-defined and hence \eqref{3.5} can then be written as
$\Xi_{1,1}(\ze,n)=(1/\pi)\Arg(M_{1,1}(\ze,n))$. Moreover, in the context of CMV operators 
we will use the convention that whenever the phrase a.e.\ is used without
further qualification, it always refers to Lebesgue measure on $\dD$.

Associated with $U$ in $\ltz$, we also introduce the two
half-lattice CMV operators $U_{\pm,n_0}$ in
$\lt{[n_0,\pm\infty)\cap\bbZ}$ by setting $\al_{n_0}=1$ which splits
the operator $U$ into a direct sum of two half-lattice operators
$U_{-,n_0-1}$ and $U_{+,n_0}$, that is,
\begin{align}
U=U_{-,n_0-1} \oplus U_{+,n_0} \, \text{ in } \,
\lt{(-\infty,n_0-1]\cap\bbZ} \oplus \lt{[n_0,\infty)\cap\bbZ}.
\end{align}
The half-lattice Weyl--Titchmarsh m-functions associated with
$U_{\pm,n_0}$ are denoted by $m_{\pm}(\cdot,n_0)$ and
$M_{\pm}(\cdot,n_0)$,
\begin{align}
m_\pm(z,n_0) &= ((U_{\pm,n_0}+zI)(U_{\pm,n_0}-zI)^{-1})(n_0,n_0),
\quad
z\in\bbC\bs\sigma(U_{\pm,n_0}), \\
M_+(z,n_0) & = m_+(z,n_0), \quad
z\in\bbC\bs\dD, \lb{3.10}\\
M_-(z,n_0) &= \frac{\Re(1+\al_{n_0}) +
i\Im(1-\al_{n_0})m_-(z,n_0-1)}{i\Im(1+\al_{n_0}) +
\Re(1-\al_{n_0})m_-(z,n_0-1)}, \quad z\in\bbC\bs\dD. \lb{3.11}
\end{align}
Then it follows that $m_\pm(\cdot,n_0)$ and $\pm M_\pm(\cdot,n_0)$
are Caratheodory functions (cf.\ \cite[Section 2]{GZ06}). Moreover,
the function $M_{1,1}(\cdot,n_0)$ is related to the m-functions
$M_\pm(\cdot,n_0)$ by (cf. \cite[Lemma 3.2]{GZ06})
\begin{equation}
M_{1,1}(z,n_0) =
\frac{1-M_+(z,n_0)M_-(z,n_0)}{M_+(z,n_0)-M_-(z,n_0)}. \lb{3.13}
\end{equation}

Next, we introduce a special class of reflectionless CMV operators
associated with a Lebesgue measurable set $\cE\subseteq\dD$ of positive Lebesgue
measure (cf. \cite{GZ06a} for a similar definition).

\begin{definition}  \lb{d3.2}
Let $\cE\subseteq\dD$ be of positive Lebesgue measure. Then we call
$U$ reflectionless on $\cE$ if for some $($equivalently, for all\,$)$ $n_0\in\bbZ$, 
\begin{align}
M_+(\ze,n_0) = -\ol{M_-(\ze,n_0)} \, \text{ for a.e.\ } \ze\in\cE.
\lb{3.19}
\end{align}
We will denote by $\cR(\cE)$ the class of all CMV operators $U$
reflectionless on $\cE$.
\end{definition}

We note that if $U$ is reflectionless on $\cE$, then by \eqref{3.5},
\eqref{3.13}, and \eqref{3.19}, one has for all $n\in\bbZ$,
\begin{align}
\Xi_{1,1}(\ze,n) = 0 \, \text{ for a.e.\ } \ze\in\cE. \lb{3.20}
\end{align}

Next, we prove the following useful result:

\begin{theorem} \lb{t3.4}
For each $n\in\bbZ$, the set
\begin{equation}
\{\ze\in\dD \st -1/2<\Xi_{1,1}(\ze,n)<1/2\} = \{\ze\in\dD \st
\Re(M_{1,1}(\ze,n))>0\} \lb{3.22}
\end{equation}
is an essential support of the absolutely continuous spectrum,
$\sigma_{\ac}(U)$, of $U$. In particular, for each $n\in\bbZ$, the
absolutely continuous spectrum coincides with the essential closure
of the set in \eqref{3.22},
\begin{equation}
\sigma_{\ac}(U)=\ol{\{\ze\in\dD \st -1/2<\Xi_{1,1}(\ze,n)<1/2\}}^e.
\lb{3.23}
\end{equation}
\end{theorem}
\begin{proof}
It follows from \eqref{3.13} that for a.e.\ $\ze\in\dD$,
\begin{align} \lb{3.23a}
&\Re(M_{1,1}(\ze,n)) \no \\
&\quad = \frac{\Re(M_+(\ze,n))(1+|M_-(\ze,n)|^2) -
\Re(M_-(\ze,n))(1+|M_+(\ze,n)|^2)}{|M_+(\ze,n)-M_-(\ze,n)|^2}.
\end{align}
Since the functions $\pm M_\pm(\cdot,n)$ are Caratheodory, that is,
$\pm\Re(M_\pm(z,n))\geq0$, $z\in\bbD$, \eqref{3.10}, \eqref{3.11},
and \eqref{3.23a} implies that up to sets of measure zero,
\begin{align}
&\{\ze\in\dD \st \Re(M_{1,1}(\ze,n))>0\} \no \\
&\quad = \{\ze\in\dD \st \Re(M_+(\ze,n))>0\} \cup \{\ze\in\dD \st
\Re(M_-(\ze,n))<0\} \no
\\
&\quad = \{\ze\in\dD \st \Re(m_+(\ze,n))>0\} \cup \{\ze\in\dD \st
\Re(m_-(\ze,n-1))>0\} \lb{3.23b}
\\
&\quad = S_+(n_0) \cup S_-(n_0). \no
\end{align}
Theorem \ref{tC.5} implies that the sets $S_\pm(n_0)$ are essential
supports of $d\mu_{+,\ac}(\cdot,n)$ and $d\mu_{-,\ac}(\cdot,n-1)$,
the absolutely continuous parts of the spectral measures of
$U_{+,n}$ and $U_{-,n-1}$, respectively. Thus, by \eqref{3.22} and
\eqref{3.23b}, $\{\ze\in\dD \st -1/2<\Xi_{1,1}(\ze,n)<1/2\}$ is an
essential support of the absolutely continuous spectrum of
$U_{+,n}\oplus U_{-,n-1}$. Since $U$ is a finite-rank perturbation
of $U_{+,n}\oplus U_{-,n-1}$ and the absolutely continuous spectrum
is invariant under finite-rank perturbations, \eqref{3.23} follows from \eqref{A.35}.
\end{proof}

With this result at hand we can give the first proof of the principal result on reflectionless CMV operators:

\begin{theorem}  \lb{t3.5}
Let $\cE\subset\dD$ be of positive Lebesgue measure and
$U\in\cR(\cE)$. Then, the absolutely continuous spectrum of $U$
contains $\ol{\cE}^e$,
\begin{equation}
\si_{\ac}(U)\supseteq\ol{\cE}^e. \lb{3.24}
\end{equation}
Moreover, the absolutely continuous spectrum of $U$ has uniform
multiplicity equal to two on $\cE$.
\end{theorem}
\begin{proof}
Since $U$ is reflectionless on $\cE$, by Definition \ref{d3.2}, one
has for each $n\in\bbZ$,
\begin{equation}
\Xi_{1,1}(\ze,n)=(1/\pi)\Im[\ln(M_{1,1}(\ze,n)]= 0 \, \text{ for
a.e.\ $\ze\in\cE$.}
\end{equation}
By \eqref{3.23}, this implies
\begin{equation}
\sigma_{\ac}(U)=\ol{\{\ze\in\dD \,|\,
-1/2<\Xi_{1,1}(\ze,n_0)<1/2\}}^e \supseteq \ol{\cE}^e
\end{equation}
for some $n_0\in\bbZ$.

Equations \eqref{3.13} and \eqref{3.19} imply
\begin{equation}
M_{1,1}(\ze,n_0)=\frac{1+|M_\pm(\ze,n_0)|^2}{\pm
2\,\Re[M_\pm(\ze,n_0)]} \, \text{ for a.e.\ $\ze\in\cE$.} \lb{3.30}
\end{equation}
Finally, combining \eqref{3.19}, \eqref{3.30}, and \eqref{C.38} then
yields that the absolutely continuous spectrum of $U$ has uniform
spectral multiplicity two on $\cE$ since
\begin{equation}
\text{for a.e.\ $\ze\in\cE$, } \; 0<\pm\Re[M_\pm(\ze,n_0)]<\infty.
\end{equation}
\end{proof}

One can also give an alternative proof of \eqref{3.24} based on the
reflectionless property of $U$ as follows (still under the
assumptions of Theorem \ref{t3.5}):

\begin{proof}[Alternative proof of \eqref{3.24}]
Fix $n_0\in\bbZ$. Since $U\in\cR(\cE)$, \eqref{3.20} yields 
\begin{equation}
\Xi_{1,1}(\ze,n_0)=\f{1}{\pi}\Im[\ln(M_{1,1}(\ze,n_0))]=0 \, \text{
for a.e.\ $\ze\in\cE$}  \lb{3.34}
\end{equation}
and hence
\begin{equation}
\Im[M_{1,1}(\ze,n_0)]=0 \, \text{ for a.e.\ $\ze\in\cE$.} \lb{3.35}
\end{equation}
If there exists a measurable set $A\subseteq \cE$ of positive
Lebesgue measure, $|A|>0$, on which $\Re(M_{1,1})$ vanishes, that
is,
\begin{equation}
\Re[M_{1,1}(\ze,n_0)]=0 \, \text{ for a.e.\ $\ze\in A$,} \lb{3.36}
\end{equation}
then \eqref{3.35} and \eqref{3.36} yield
\begin{equation}
M_{1,1}(\ze,n_0)=0 \, \text{ for a.e.\ $\ze\in A$.}
\end{equation}
Since $M_{1,1}(\cdot,n_0)$ is a Caratheodory function, Theorem
\ref{tC.2}\,$(ii)$ yields the contradiction $M_{1,1}\equiv 0$. Thus,
no such set $A\subseteq\cE$ exists and one concludes that
\begin{equation}
\Re[M_{1,1}(\ze,n_0)]>0 \, \text{ for a.e.\ $\ze\in \cE$.} \lb{3.38}
\end{equation}
Moreover, since $M_{1,1}(\ze,n_0)$ exists finitely for a.e.\
$\ze\in\dD$ by Theorem \ref{tC.2}\,$(i)$, this yields
\begin{equation}
0<\Re[M_{1,1}(\ze,n_0)]<\infty \, \text{ for a.e.\ $\ze\in \cE$.}
\lb{3.39}
\end{equation}
Let $d\Om_{1,1}(\cdot,n_0)$ denote the measure in the Herglotz
representation \eqref{C.3} of the function $M_{1,1}(\cdot,n_0)$.
Then by \eqref{C.11}
\begin{align}
S_{\ac}(n_0)=\{\ze\in\dD \st 0<\Re[M_{1,1}(\ze,n_0)]<\infty\}
\end{align}
is an essential support of the absolutely continuous part
$d\Om_{1,1, \rm ac}(\cdot,n_0)$ of the measure
$d\Om_{1,1}(\cdot,n_0)$. Thus \eqref{3.39} implies
\begin{equation}
\supp\,[d\Om_{1,1, \rm ac}(\cdot,n_0)] = \ol{S_{\ac}(n_0)}^e
\supseteq\ol{\cE}^e. \lb{3.42}
\end{equation}

Next, we introduce the $2\times 2$ matrix-valued Weyl--Titchmarsh
$m$-function associated with $U$,
\begin{align}
& M(z,n_0) =
\begin{pmatrix}
M_{0,0}(z,n_0) & M_{0,1}(z,n_0) \\
M_{1,0}(z,n_0) & M_{1,1}(z,n_0)
\end{pmatrix} \no
\\
&\quad = \begin{pmatrix} (\de_{n_0-1},(U+zI)(U-zI)^{-1}\de_{n_0-1})
&(\de_{n_0-1},(U+zI)(U-zI)^{-1}\de_{n_0})
\\
(\de_{n_0},(U+zI)(U-zI)^{-1}\de_{n_0-1}) &
(\de_{n_0},(U+zI)(U-zI)^{-1}\de_{n_0})
\end{pmatrix} \no
\\
&\quad = \oint_\dD d\Omega(\ze,n_0)\, \frac{\ze+z}{\ze-z}, \quad
z\in\bbD, \lb{3.42a}
\end{align}
where $d\Om=(d\Om_{j,k})_{j,k=0,1}$ denotes a $2\times 2$
matrix-valued nonnegative measure satisfying
\begin{equation}
\oint_\dD d\,|\Om_{j,k}(\ze)| < \infty, \quad j,k=0,1.
\end{equation}
It is proven in \cite[Corollary 3.5]{GZ06} that $d\Om(\cdot,n_0)$ is
the spectral measure of $U$, that is, $U$ is unitarily equivalent to
the operator of multiplication by $I_2\id$ $($where $I_2$ is the
$2\times 2$ identity matrix and $\id(\ze)=\ze$, $\ze\in\dD$$)$ on
$L^2(\dD; d\Om(\cdot,n_0))$. Then since the matrix-valued measure
$d\Om(\cdot,n_0)$ and its trace measure $d\Om^{\tr}(\cdot,n_0) =
d\Om_{0,0}(\cdot,n_0) + d\Om_{1,1}(\cdot,n_0)$ are mutually
absolutely continuous, one has
\begin{align}
\si_\ac(U) = \supp(d\Om_\ac(\cdot,n_0)) =
\supp(d\Om^\tr_\ac(\cdot,n_0)). \lb{3.42b}
\end{align}
Note that by \eqref{3.42a}
\begin{align}
M_{0,0}(z,n_0) = M_{1,1}(z,n_0-1),
\end{align}
and hence
\begin{align}
d\Om^{\tr}(\cdot,n_0) = d\Om_{1,1}(\cdot,n_0-1) +
d\Om_{1,1}(\cdot,n_0). \lb{3.42c}
\end{align}
Thus, it follows from \eqref{3.42}, \eqref{3.42b}, and \eqref{3.42c}
that
\begin{equation}
\sigma_{\ac}(U) =
\supp\,\big[d\Om^{\tr}_{\ac}(\cdot,n_0)\big]\supseteq\ol{\cE}^e.
\lb{3.43}
\end{equation}
\end{proof}

We note that the strategy of proof in \eqref{3.34}--\eqref{3.43} is well-known. It goes back to Kotani \cite{Ko84}--\cite{Ko87b} and has been exploited in \cite[Theorem 12.5]{PF92}.

\medskip

Next, we briefly turn to Jacobi operators on $\bbZ$. 

We start with some general considerations of self-adjoint Jacobi
operators. Let $a=\{a(n)\}_{n\in\bbZ}$ and $b=\{b(n)\}_{n\in\bbZ}$ be
two sequences (Jacobi parameters) satisfying
\begin{equation}
a, b \in \ell^\infty(\bbZ), \quad a(n)>0, \; b(n)\in\bbR, \; n\in\bbZ,
\end{equation}
and denote by $L$ the second-order difference expression defined by
\begin{equation}
L = a S^+ + a^- S^- + b.
\end{equation}
Moreover, we introduce the associated bounded self-adjoint Jacobi
operator $H$ in $\ltz$ by
\begin{align}
&(Hf)(n) = (L f)(n), \quad n\in\bbZ, \no
\\
&f=\{f(n)\}_{n\in\bbZ}\in\dom(H)=\ltz. \lb{2.2}
\end{align}
Next, let $g(z,\cdot)$ denote the diagonal Green's function of
$H$, that is,
\begin{align}
g(z,n)=G(z,n,n), \quad G(z,n,n')=(H-zI)^{-1}(n,n'),& \no
\\
z\in\bbC\bs\sigma(H), \; n,n'\in\bbZ.&
\end{align}
Since for each $n\in\bbZ$, $g(\cdot,n)$ is a Herglotz function
(i.e., it maps the open complex upper half-plane analytically to
itself),
\begin{equation}
\xi(\la,n)=\f{1}{\pi}\lim_{\eps\downarrow 0}
\Im[\ln(g(\la+i\eps,n))] \, \text{ for a.e.\ $\la\in\bbR$} \lb{2.5}
\end{equation}
is well-defined for each $n\in\bbZ$. In particular, for all
$n\in\bbZ$,
\begin{equation}
0 \leq \xi(\la,n) \leq 1 \, \text{ for a.e.\ $\la\in\bbR$.}
\end{equation}

In the following we will frequently use the convenient abbreviation
\begin{equation}
h(\la_0+i0)=\lim_{\eps\downarrow 0} h(\la_0 +i\eps), \quad
\la_0\in\bbR, \lb{2.7}
\end{equation}
whenever the limit in \eqref{2.7} is well-defined and hence
\eqref{2.5} can then be written as
$\xi(\la,n)=(1/\pi)\Arg(g(\la+i0,n))$. Moreover, for the remainder of this section we
will use the convention that whenever the phrase a.e.\ is used
without further qualification, it always refers to Lebesgue measure
on $\bbR$.

Associated with $H$ in $\ltz$, we also introduce the two half-lattice
Jacobi operators $H_{\pm,n_0}$ in $\lt{[n_0,\pm\infty)\cap\bbZ}$ by
\begin{align}
&H_{\pm,n_0} = P_{\pm,n_0} H
P_{\pm,n_0}|_{\lt{[n_0,\pm\infty)\cap\bbZ}},
\end{align}
where $P_{\pm,n_0}$ are the orthogonal projections onto the subspaces
$\lt{[n_0,\pm\infty)\cap\bbZ}$. By inspection, $H_{\pm,n_0}$ satisfy
Dirchlet boundary conditions at $n_0\mp1$, that is,
\begin{align}
&(H_{\pm,n_0}f)(n) = (Lf)(n), \quad n\gtreqless n_0, \no
\\
&f\in\dom(H_{\pm,n_0})=\lt{[n_0,\pm\infty)\cap\bbZ}, \quad
f(n_0\mp1)=0.
\end{align}

The half-lattice Weyl--Titchmarsh m-functions associated with
$H_{\pm,n_0}$ are denoted by $m_{\pm}(\cdot,n_0)$ and
$M_{\pm}(\cdot,n_0)$,
\begin{align}
m_\pm(z,n_0) &= (\delta_{n_0}, (H_{\pm,n_0}-zI)^{-1}
\delta_{n_0})_{\lt{[n_0,\pm\infty)\cap\bbZ}}, \quad
z\in\bbC\bs\sigma(H_{\pm,n_0}),   \lb {2.9} \\
M_+(z,n_0) & = -m_+(z,n_0)^{-1}-z+b(n_0), \quad
z\in\bbC\bs\bbR, \lb{2.10}\\
M_-(z,n_0) & = m_-(z,n_0)^{-1}, \quad z\in\bbC\bs\bbR, \lb{2.11}
\end{align}
where $\delta_k = \{\delta_{k,n}\}_{n\in\bbZ}$, $k\in\bbZ$.
An equivalent definition of $M_\pm(\cdot,n_0)$ is
\begin{align}
M_\pm(z,n_0)=-a(n_0)\frac{\psi_\pm(z,n_0+1)}{\psi_\pm(z,n_0)}, \quad
z\in\bbC\bs\bbR, \lb{2.12}
\end{align}
where $\psi_\pm(z,\cdot)$ are the Weyl--Titchmarsh solutions of
$(L-z)\psi_\pm(z,\cdot)=0$ with
$\psi_\pm(z,\cdot)\in\lt{[n_0,\pm\infty)\cap\bbZ}$. Then it follows
that the diagonal Green's function $g(\cdot,n_0)$ is related to the
m-functions $M_\pm(\cdot,n_0)$ via
\begin{align}
g(z,n_0) &= [M_-(z,n_0)-M_+(z,n_0)]^{-1}.   \lb{2.13}
\end{align}

Next, following \cite{GKT96}, we introduce a special class of
reflectionless Jacobi operators associated with a Lebesgue
measurable set $\cE\subset\bbR$ of positive Lebesgue measure (cf. also
\cite{Re07} and \cite[Lemma 8.1]{Te00}).

\begin{definition}  \lb{d2.2}
Let $\cE\subset\bbR$ be of positive Lebesgue measure. Then we call
$H$ reflectionless on $\cE$ if for some $($equivalently, for all\,$)$ $n_0\in\bbZ$, 
\begin{align}
M_+(\la+i0,n_0) = \ol{M_-(\la+i0,n_0)} \, \text{ for a.e.\ }
\la\in\cE. \lb{2.19}
\end{align}
Equivalently, $H$ is called reflectionless on $\cE$ if for all
$n\in\bbZ$,
\begin{align}
\xi(\la,n) = 1/2 \, \text{ for a.e.\ } \la\in\cE. \lb{2.20}
\end{align}
We will denote by $\cR(\cE)$ the class of all Jacobi operators $H$
reflectionless on $\cE$.
\end{definition}

We recall the following result on essential supports of
the absolutely continuous spectrum of Jacobi operators originally proven in
\cite{GS96} (see also \cite[Lemma 3.11, (B.28)]{Te00}):

\begin{theorem}[\cite{GS96}] \lb{t2.4}
For each $n\in\bbZ$, the set
\begin{equation}
\{\la\in\bbR\,|\, 0<\xi(\la,n)<1\}  \lb{2.22}
\end{equation}
is an essential support of the absolutely continuous spectrum,
$\sigma_{\ac}(H)$, of $H$. In particular, for each $n\in\bbZ$, the
absolutely continuous spectrum coincides with the essential closure
of the set in \eqref{2.22},
\begin{equation}
\sigma_{\ac}(H)=\ol{\{\la\in\bbR\,|\, 0<\xi(\la,n)<1\}}^e. \lb{2.23}
\end{equation}
\end{theorem}

With this result at hand we can give a short proof of the
principal result on reflectionless Jacobi operators:

\begin{theorem}  \lb{t2.5}
Let $\cE\subset\bbR$ be of positive Lebesgue measure and
$H\in\cR(\cE)$. Then, the absolutely continuous spectrum of $H$
contains $\ol{\cE}^e$,
\begin{equation}
\si_{\ac}(H)\supseteq\ol{\cE}^e. \lb{2.24}
\end{equation}
Moreover, the absolutely continuous spectrum of $H$ has uniform
multiplicity equal to two on $\cE$.
\end{theorem}
\begin{proof}
Since $H$ is reflectionless on $\cE$, by Definition \ref{d2.2}, one
has for each $n\in\bbZ$,
\begin{equation}
\xi(\la,n)=(1/\pi)\Im[\ln(g(\la+i0,n)]= 1/2 \, \text{ for a.e.\
$\la\in\cE$.}
\end{equation}
By \eqref{2.23} and Lemma \ref{lA.7}, this implies
\begin{equation}
\sigma_{\ac}(H)=\ol{\{\la\in\bbR \,|\, 0<\xi(\la,n_0)<1\}}^e
\supseteq \ol{\cE}^e
\end{equation}
for some $n_0\in\bbZ$.

Equations \eqref{2.13} and \eqref{2.19} imply
\begin{equation}
-1/g(\la+i0,n_0)=\pm 2i \, \Im[M_\pm(\la+i0,n_0)] \, \text{ for
a.e.\ $\la\in\cE$.}    \lb{2.30}
\end{equation}
Finally, combining \eqref{2.19}, \eqref{2.30}, and \eqref{B.38} then
yields that the absolutely continuous spectrum of $H$ has uniform
spectral multiplicity two on $\cE$ since
\begin{equation}
\text{for a.e.\ $\la\in\cE$, } \;
0<\pm\Im[M_\pm(\la+i0,n_0)]<\infty.
\end{equation}
\end{proof}

Next, we briefly turn to Schr\"odinger operators on $\bbR$. Let
\begin{equation}
V\in L^\infty(\bbR;dx), \quad V \, \text{ real-valued,}
\end{equation}
and consider the differential expression
\begin{equation}
L=-d^2/dx^2+V(x), \quad x\in\bbR.
\end{equation}
We denote by $H$ the corresponding self-adjoint
realization of $L$ in $L^2(\bbR;dx)$ given by
\begin{equation}
Hf=L f, \quad f\in\dom(H)=H^2(\bbR), 
\end{equation}
with $H^2(\bbR)$ the usual Sobolev space. Let $g(z,\cdot)$ denote the diagonal Green's function of $H$, that is,
\begin{equation}
g(z,x)=G(z,x,x), \quad G(z,x,x')=(H-zI)^{-1}(x,x'), \quad
z\in\bbC\backslash\sigma(H), \; x,x'\in\bbR.
\end{equation}
Since for each $x\in\bbR$, $g(\cdot,x)$ is a Herglotz function,
\begin{equation}
\xi(\lambda,x)=\f{1}{\pi}\lim_{\varepsilon\downarrow 0}
\Im[\ln(g(\lambda+i\varepsilon,x))] \, \text{ for a.e.\
$\lambda\in\bbR$}  \lb{4.5}
\end{equation}
is well-defined for each $x\in\bbR$. In particular, for all $x\in\bbR$,
\begin{equation}
0 \leq \xi(\lambda,x) \leq 1 \, \text{ for a.e.\ $\lambda\in\bbR$.}
\end{equation}

Associated with $H$ in $L^2(\bbR;dx)$ we also introduce the two half-line
Schr\"odinger operators $H_{\pm,x_0}$ in $L^2([x_0,\pm\infty);dx)$ with
Dirchlet boundary conditions at the finite endpoint $x_0\in\bbR$,
\begin{align}
& H_{\pm,x_0} f= Lf, \no \\
& f\in\dom(H_{\pm,x_0})=\big\{g\in L^2([x_0,\pm\infty);dx)\,|\, g, g'\in
AC([x_0, x_0\pm R]) \, \text{for all $R>0$;} \no \\
& \hspace*{4.5cm} \lim_{\varepsilon\downarrow 0}
g(x_0\pm\varepsilon)=0; \, Lg\in L^2([x_0,\pm\infty);dx)\big\}.
\end{align}
Denoting by $\psi_{\pm}(z,\cdot)$ the Weyl--Titchmarsh solutions of
$(L - z)\psi(z,\dott) = 0$, satisfying
\begin{equation}
\psi_{\pm}(z,\cdot)\in L^2([x_0,\pm\infty);dx),
\end{equation}
the half-line Weyl--Titchmarsh functions associated with
$H_{\pm,x_0}$ are given by
\begin{equation}
m_{\pm}(z,x_0)=\f{\psi_{\pm}'(z,x_0)}{\psi_{\pm}(z,x_0)}, \quad
z\in\bbC\backslash\sigma(H_{\pm,x_0}). \lb{4.10}
\end{equation}
Then the diagonal Green's function of $H$ satisfies
\begin{equation}
g(z,x_0)=[m_{-}(z,x_0) - m_{+}(z,x_0)]^{-1}.    \lb{4.11}
\end{equation}

We note that the condition $V\in L^\infty(\bbR)$ in this section is only used 
for simplicity. The general case $V \in L^1_{\loc}(\bbR)$ and $L$ in the 
limit point case at $\pm\infty$ is discussed in detail in \cite{GY06}. 

Next, we introduce a special class of reflectionless Schr\"odinger
operators associated with a Lebesgue measurable set $\cE\subset\bbR$ of
positive Lebesgue measure.

\begin{definition}  \lb{d4.2}
Let $\cE\subset\bbR$ be of positive Lebesgue measure and pick
$x_0\in\bbR$. Then $H$ is called reflectionless on $\cE$ if
\begin{equation}
m_+(\lambda+i0,x_0)= \ol{m_-(\lambda+i0,x_0)} \, \text{ for a.e.\
$\lambda\in\cE$.}  \lb{4.19}
\end{equation}
Equivalently, $H$ is called reflectionless if for each $x\in\bbR$,
\begin{equation}
 \xi(\lambda,x)=1/2 \, \text{ for a.e.\ $\lambda\in\cE$.}    \lb{4.20}
\end{equation}
We will denote by $\cR(\cE)$ the class of Schr\"odinger operators
$H$ reflectionless on $\cE$.
\end{definition}

We recall the following result on essential supports of
the absolutely continuous spectrum of Schr\"odinger operators on the
real line proven in \cite{GS96} (see also  \cite[p.\ 383]{AD56}): 

\begin{theorem} [\cite{GS96}] \lb{t4.4}
For each $x\in\bbR$, the set
\begin{equation}
\{\lambda\in\bbR\,|\, 0<\xi(\lambda,x)<1\}
\end{equation}
is an essential support of the absolutely continuous spectrum,
$\sigma_{\rm ac}(H)$, of $H$. In particular, for each $x\in\bbR$,
\begin{equation}
\sigma_{\rm ac}(H)=\ol{\{\lambda\in\bbR\,|\,
0<\xi(\lambda,x)<1\}}^e. \lb{4.23}
\end{equation}
\end{theorem}

Given Theorem \ref{t4.4}, the principal result for reflectionless Schr\"odinger operators  then reads as follows: 

\begin{theorem}  \lb{t4.5}
Let $\cE\subset\bbR$ be of positive Lebesgue measure and
$H\in\cR(\cE)$. Then, the absolutely continuous spectrum of $H$
contains $\ol{\cE}^e$,
\begin{equation}
\sigma_{\rm ac}(H) \supseteq \ol{\cE}^e.   \lb{4.24}
\end{equation}
Moreover, the absolutely continuous spectrum of $H$ has uniform
multiplicity equal to two on $\cE$.
\end{theorem}
\begin{proof}
Since $H$ is reflectionless, one has for each $x\in\bbR$,
\begin{equation}
\xi(\lambda,x)=(1/\pi)\Im[\ln(g(\lambda+i0,x)]= 1/2 \, \text{ for
a.e.\ $\lambda\in\cE$.}
\end{equation}
By \eqref{4.23} and Lemma \ref{lA.7}, this implies
\begin{equation}
\sigma_{\rm ac}(H)=\ol{\{\lambda\in\bbR \,|\,
0<\xi(\lambda,x_0)<1\}}^e \supseteq \ol{\cE}^e
\end{equation}
for some $x_0\in\bbR$.

Equations \eqref{4.11} and \eqref{4.19} imply
\begin{equation}
-1/g(\lambda+i0,x_0)=\pm 2i \, \Im[m_\pm(\lambda+i0,x_0)] \, \text{
for a.e.\ $\lambda\in\cE$.}    \lb{4.30}
\end{equation}
Finally, combining \eqref{4.19}, \eqref{4.30}, and \eqref{B.38} then
yields that the absolutely continuous spectrum of $H$ has uniform
spectral multiplicity two on $\cE$ since
\begin{equation}
\text{for a.e.\ $\lambda\in\cE$, } \;
0<\pm\Im[m_\pm(\lambda+i0,x_0)]<\infty.
\end{equation}
\end{proof}

\begin{remark} \lb{r4.6}   
$(i)$ As in the case of CMV operators one can formulate alternative proofs of 
Theorems \ref{t2.5} and \ref{t4.5} based on the reflectionless property of $H$ precisely 
along the lines of \eqref{3.34}--\eqref{3.43}. \\
$(ii)$ One particularly interesting situation occurs in connection with Definition \ref{d3.2} when $\si(U)=\si_{\ess}(U)$ is a homogeneous set (cf.\ \cite{Ca83},
\cite{SY95}, \cite{SY97} for the definition of homogeneous sets) and $U$ is
reflectionless on $\si(U)$. This case has been studied by Peherstorfer and Yuditskii 
\cite{PY03}, \cite{PY06}, and more recently, in \cite{GZ08}. The same applies to Definitions \ref{d2.2} and \ref{d4.2} in the case of Jacobi and Schr\"odinger operators which were studied by Sodin and Yuditskii \cite{SY95a}--\cite{SY96}, and more
recently, in \cite{GY06}, \cite{GZ08}. \\
$(iii)$ As shown in \cite{GY06} for reflectionless Schr\"odinger operators with 
$\sigma(H)=\cE$, and in \cite{GZ08} for reflectionless CMV, Jacobi, and Schr\"odinger operators on $\cE$, under some additional assumptions on $\cE$ (such as $\cE$ a homogeneous set, etc.), it is possible to prove the absence of any singular spectrum of $U$ and $H$ on $\cE$. But unlike the most elementary proofs presented in this
section, the results in \cite{GY06} and \cite{GZ08} rely on sophisticated
techniques due to Zinsmeister \cite{Zi89} and Peherstorfer and Yuditskii \cite{PY03}, 
\cite{PY06}, respectively.
\end{remark}

\appendix
\section{Herglotz Functions and Weyl--Titchmarsh Theory
for Jacobi and Schr\"odinger Operators in a Nutshell} \lb{sB}
\renewcommand{\theequation}{A.\arabic{equation}}
\renewcommand{\thetheorem}{A.\arabic{theorem}}
\setcounter{theorem}{0}
\setcounter{equation}{0}

The material in this appendix is known, but since we use it repeatedly at
various places in Section \ref{s3}, we thought it worthwhile to collect it
in an appendix.

\begin{definition} \lb{dB.1}
Let $\bbC_\pm=\{z\in\bbC \mid \Im(z)\gtrless 0 \}$.
$m:{\mathbb{C_+}}\to {\mathbb{C}}$ is called a {\it Herglotz}
function $($or {\it Nevanlinna} or {\it Pick} function$)$ if $m$ is analytic
on
${\mathbb{C}}_+$ and
$m({\mathbb{C}}_+)\subseteq {\mathbb{C}}_+$.
\end{definition}

\smallskip

\noindent One then extends $m$ to $\bbC_-$ by reflection, that is, one
defines
\begin{equation}
m(z)=\overline{m(\overline z)},
\quad z\in{\mathbb{C}}_-. \lb{B.1}
\end{equation}
Of course, generally, \eqref{B.1} does not represent the analytic
continuation of $m\big|_{\bbC_+}$ into $\bbC_-$.

The fundamental result on Herglotz functions and their representations
on Borel transforms, in part due to Fatou, Herglotz, Luzin,
Nevanlinna, Plessner, Privalov, de la Vall{\'e}e Poussin, Riesz, and
others, then reads as follows.

\begin{theorem} $($\cite[Sect.~69]{AG81}, \cite{AD56},
\cite[Chs.~II, IV]{Do74}, \cite{KK74}, \cite[Ch.~6]{Ko98},
\cite[Chs.~II, IV]{Pr56}, \cite[Ch.~5]{RR94}$)$. \lb{tB.2}
Let $m$ be a Herglotz function. Then, \\
$(i)$ $m(z)$ has finite normal limits $m(\la \pm i0)=\lim_{\eps
\downarrow 0} m(\la \pm i\eps)$ for
a.e.~$\la \in {\mathbb{R}}$. \\
$(ii)$ Suppose $m(z)$ has a zero normal limit on a subset of
${\mathbb{R}}$ having positive Lebesgue measure. Then $m\equiv 0$. \\
$(iii)$ There exists a nonnegative measure $d\om$ on ${\mathbb{R}}$
satisfying
\begin{equation} \lb{B.2}
\int_{{\mathbb{R}}} \frac{d\om (\la )}{1+\la^2}<\infty
\end{equation}
such that the Nevanlinna, respectively, Riesz--Herglotz
representation
\begin{align}
\begin{split}
&m(z)=c+dz+\int_{{\mathbb{R}}} d\om (\la) \bigg(\frac{1}{\la
-z}-\frac{\la}
{1+\la^2}\bigg), \quad z\in\bbC_+, \lb{B.3} \\[2mm]
& \, c=\Re[m(i)],\quad d=\lim_{\eta \uparrow
\infty}m(i\eta )/(i\eta ) \geq 0
\end{split}
\end{align}
holds. Conversely, any function $m$ of the type \eqref{B.3} is a
Herglotz function. \\
$(iv)$ The absolutely continuous $({\it ac})$ part $d\om_{ac}$ of
$d\om$ with respect to Lebesgue measure $d\la$ on ${\mathbb{R}}$ is
given by
\begin{equation}\lb{B.5}
d\om_{ac}(\la)=\pi^{-1}\Im[m(\la+i0)]\,d\la.
\end{equation}
\end{theorem}

Next, we denote by
\begin{equation}
d\om =d\om_{\ac}+d\om_{\sc} +d\om_{\pp} \lb{B.8}
\end{equation}
the decomposition of $d\om$ into its absolutely continuous $({\it
ac})$, singularly continuous $({\it sc})$, and pure point $({\it
pp})$ parts with respect to Lebesgue measure on $\bbR$.

\begin{theorem} $($\cite{Gi84}, \cite{GP87}$)$.  \lb{tB.5}
Let $m$ be a Herglotz function with representation \eqref{B.3} and
denote by $\La$ the set
\begin{equation}
\La=\{\la\in\bbR\,|\, \Im[m(\la+i0)] \, \text{exists $($finitely or
infinitely$)$}\}. \lb{B.9}
\end{equation}
Then, $S$, $S_{\ac}$, $S_{\rm s}$, $S_{\sc}$, $S_{\pp}$ are
essential supports of $d\om$, $d\om_{\ac}$, $d\om_{\rm s}$,
$d\om_{\sc}$, $d\om_{\pp}$, respectively, where
\begin{align}
S&=\{\la\in\La\,|\, 0<\Im[m(\la+i0)]\leq\infty\}, \lb{B.10}
\\
S_{\ac}&=\{\la\in\La\,|\, 0<\Im[m(\la+i0)]<\infty\},
\lb{B.11} \\
S_{\rm s}&=\{\la\in\La\,|\, \Im[m(\la+i0)]=\infty\},
\lb{B.12} \\
S_{\sc}&=\Big\{\la\in\La \,\Big|\, |\, \Im[m(\la+i0)]=\infty, \,
\lim_{\eps\downarrow 0}
(-i\eps)m(\la+i\eps)=0\Big\}, \lb{B.13} \\
S_{\pp}&=\Big\{\la\in\La \,\Big|\, \Im[m(\la+i0)]=\infty, \,
\lim_{\eps\downarrow 0} (-i\eps)m(\la+i\eps)=\om(\{\la\})>0\Big\}.
\lb{B.14}
\end{align}
Moreover, since
\begin{align}
& |\{\la\in\bbR\,|\, |m(\la+i0)|=\infty\}|=0, \lb{B.15}  \\
& \om(\{\la\in\bbR\,|\, |m(\la+i0)|=\infty, \;
\Im[m(\la+i0)]<\infty\})=0,  \lb{B.16}
\end{align}
also $S_{\rm s}^{\prime}$, $S_{\sc}^{\prime}$, $S_{\pp}^{\prime}$
are essential supports of $d\om_{\rm s}$, $d\om_{\sc}$,
$d\om_{\pp}$, respectively, where
\begin{align}
S_{\rm ac}^{\prime}&=\{\la\in\La\,|\, \Im[m(\la+i0)]>0\},
\lb{B.17a} \\
S_{\rm s}^{\prime}&=\{\la\in\La\,|\, |m(\la+i0)|=\infty\},
\lb{B.17} \\
S_{\sc}^{\prime}&=\Big\{\la\in\La \,\Big|\, |m(\la+i0)|=\infty, \,
\lim_{\eps\downarrow 0}
(-i\eps)m(\la+i\eps)=0\Big\}, \lb{B.18} \\
S_{\pp}^{\prime}&=\Big\{\la\in\La \,\Big|\, |m(\la+i0)|=\infty, \,
\lim_{\eps\downarrow 0} (-i\eps)m(\la+i\eps)=\om(\{\la\})>0\Big\}.
\lb{B.19}
\end{align}
In particular $($cf.\ Lemma \ref{lA.3}$)$,
\begin{equation}
S_{\rm s}\sim S_{\rm s}^{\prime}, \quad S_{\sc}\sim S_{\rm
sc}^{\prime}, \quad S_{\pp}\sim S_{\pp}^{\prime}. \lb{B.20}
\end{equation}
\end{theorem}

Next, consider Herglotz functions $\pm m_\pm$ of the type \eqref{B.3},
\begin{align}
\begin{split}
& \pm m_\pm(z)=c_\pm+d_\pm z+\int_{{\mathbb{R}}} d\om_\pm (\la)
\bigg(\frac{1}{\la -z}-\frac{\la}
{1+\la^2}\bigg), \quad z\in\bbC_+, \lb{B.21} \\
& \,\, c_\pm\in\bbR, \quad d_\pm\geq 0,
\end{split}
\end{align}
and introduce the $2\times 2$ matrix-valued Herglotz function $M$
\begin{align}
&M(z)=\big(M_{j,k}(z)\big)_{j,k=0,1}, \quad z\in\bbC_+, \lb{B.22} \\
&M(z)=\f{1}{m_-(z)-m_+(z)}\begin{pmatrix} 1 & \f{1}{2}[m_-(z)+m_+(z)]
\\
\f{1}{2}[m_-(z)+m_+(z)] & m_-(z)m_+(z) \end{pmatrix}
\lb{B.23} \\
& \hspace*{.85cm} =C+Dz + \int_{{\mathbb{R}}} d\Om (\la)
\bigg(\frac{1}{\la -z}-\frac{\la}
{1+\la^2}\bigg), \quad z\in\bbC_+,  \lb{B.24} \\
& \, C=C^*, \quad D \geq 0 \no
\end{align}
with $C=(C_{j,k})_{j,k=0,1}$ and $D=(D_{j,k})_{j,k=0,1}$ $2\times 2$
matrices and $d\Om=(d\Om_{j,k})_{j,k=0,1}$ a $2\times 2$
matrix-valued nonnegative measure satisfying
\begin{equation}
\int_\bbR \f{d|\Om_{j,k}(\la)|}{1+\la^2} < \infty, \quad j,k=0,1.
\lb{B.25}
\end{equation}
Moreover, we introduce the trace Herglotz function $M^{\tr}$
\begin{align}
&M^{\tr}(z)=M_{0,0}(z)+M_{1,1}(z)=\f{1+m_-(z)m_+(z)}{m_-(z)-m_+(z)}
\lb{B.26} \\
& \hspace*{1.05cm} =c + dz + \int_{\bbR} d\Om^{\tr}(\la)
\bigg(\frac{1}{\la -z}-\frac{\la} {1+\la^2}\bigg), \quad z\in\bbC_+,
\lb{B.27} \\ & \, c\in\bbR, \;  d\geq 0, \quad d\Om^{\tr}=
d\Om_{0,0} + d\Om_{1,1}. \no
\end{align}
Then (cf., e.g., \cite[p.\ 21]{CL90}),
\begin{equation}
d\Om \ll d\Om^{\tr} \ll d\Om   \lb{B.28}
\end{equation}
(where $d\mu \ll d\nu$ denotes that $d\mu$ is absolutely
continuous with respect to $d\nu$). \eqref{B.28} follows from the 
fact that for any nonnegative $2\times 2$ matrix 
$A=\big(A_{j,k}\big)_{1\leq j,k\leq2}$ with complex-valued 
entries, 
\begin{equation}
|A_{j,k}| \leq A_{j,j}^{1/2} A_{k,k}^{1/2} \leq (A_{j,j} + A_{k,k})/2, \quad 
1 \leq j,k \leq 2,
\end{equation}
and hence
\begin{equation}
d \Omega_{j,k} \ll d \Omega_{j,j} + d \Omega_{k,k} \ll d\Om^{\tr}, \quad 
1 \leq j,k \leq 2. 
\end{equation}

The next result holds for the Jacobi and Schr\"odinger cases.
In the Jacobi case we identify
\begin{equation}
m_\pm(z) \, \text{ and } \, M_\pm(z,n_0), \quad z\in\bbC_+,
\end{equation}
where $M_\pm(z,n_0)$ denote the half-lattice Weyl--Titchmarsh
$m$-functions defined in \eqref{2.10}--\eqref{2.12} and in the
Schr\"odinger case we identify 
\begin{equation}
m_\pm(z) \, \text{ and } \, m_\pm(z,x_0), \quad z\in\bbC_+,
\end{equation}
where $m_\pm(z,x_0)$ are the half-line Weyl--Titchmarsh
$m$-functions defined in \eqref{4.10}. One then has the following
basic result.

\begin{theorem} $($\cite{Gi98}, \cite{Ka62}, \cite{Ka63}, \cite{Si05b}$)$. 
\lb{tB.7} ${}$ \\
$(i)$ The spectral multiplicity of the Jacobi or Schr\"odinger operator $H$ 
is two if and only if
\begin{equation}
|\cM_2|>0,  \lb{B.37}
\end{equation}
where
\begin{equation}
\cM_2=\{\la\in\La_+\,|\,
m_+(\la+i0)\in\bbC\bs\bbR\}\cap\{\la\in\La_-\,|\,
m_-(\la+i0)\in\bbC\bs\bbR\}.  \lb{B.38}
\end{equation}
If $|\cM_2|=0$, the spectrum of $H$ is simple. Moreover, $\cM_2$ is
a maximal set on which $H$ has uniform multiplicity two. \\
$(ii)$ A maximal set $\cM_1$ on which $H$ has uniform multiplicity
one is given by
\begin{align}
\cM_1&=\{\la\in\La_+\cap\La_-\,|\, m_+(\la+i0)=
m_-(\la+i0)\in\bbR\}  \no \\
&\quad \cup \{\la\in\La_+\cap\La_-\,|\,
|m_+(\la+i0)|= |m_-(\la+i0)|=\infty\}  \no \\
&\quad \cup \{\la\in\La_+\cap\La_-\,|\, m_+(\la+i0)\in\bbR,
m_-(\la+i0)\in\bbC\bs\bbR\}  \no \\
&\quad \cup \{\la\in\La_+\cap\La_-\,|\, m_-(\la+i0)\in\bbR,
m_+(\la+i0)\in\bbC\bs\bbR\}. \lb{B.39}
\end{align}
In particular, $\sigma_{\rm s}(H)=\sigma_{\sc}(H)\cup
\sigma_{\pp}(H)$ is always simple.
\end{theorem}

\section{Caratheodory Functions and Weyl--Titchmarsh Theory
for CMV Operators in a Nutshell} \lb{sC}
\renewcommand{\theequation}{B.\arabic{equation}}
\renewcommand{\thetheorem}{B.\arabic{theorem}}
\setcounter{theorem}{0} \setcounter{equation}{0}

In this appendix we provide some basic facts on Caratheodory
functions and prove the analog of Theorem \ref{tB.7} for CMV
operators.

\begin{definition} \lb{dC.1}
Let $\bbD$ and $\dD$ denote the open unit disk and the
counterclockwise oriented unit circle in the complex plane $\bbC$,
\begin{equation}
\bbD = \{ z\in\bbC \st \abs{z} < 1 \}, \quad \dD = \{ \ze\in\bbC \st
\abs{\ze} = 1 \},
\end{equation}
and $\Cl$ and $\Cr$ the open left and right complex half-planes,
respectively,
\begin{equation}
\Cl = \{z\in\bbC \st \Re(z) < 0\}, \quad \Cr = \{z\in\bbC \st \Re(z)
> 0\}.
\end{equation}
A function $f:\bbD\to\bbC$ is called Caratheodory if $f$ is analytic
on $\bbD$ and $f(\bbD)\subset\Cr$. One then extends $f$ to
$\bbC\bs\ol{\bbD}$ by reflection, that is, one defines
\begin{align}
f(z) = -\ol{f(1/\ol{z})}, \quad z\in\bbC\bs\ol{\bbD}. \lb{C.1}
\end{align}
Of course, generally, \eqref{C.1} does not represent the analytic
continuation of $f\big|_\bbD$ into $\bbC\bs\ol{\bbD}$.
\end{definition}

The fundamental result on Caratheodory functions then reads as follows:

\begin{theorem} $($\cite[Sect.\ 3.1]{Ak65}, \cite[Sect.\ 69]{AG81},
\cite[Sect.\ 1.3]{Si05}$)$. \lb{tC.2} Let $f$ be a Caratheodory function. Then, \\
$(i)$ $f(z)$ has finite normal limits $f(\ze)=\lim_{r\uparrow1}
f(r\ze)$ for a.e.~$\ze\in\dD$. \\
$(ii)$ Suppose $f(r\ze)$ has a zero normal limit on a subset of
$\dD$ having positive Lebesgue measure as $r\uparrow 1$. Then $f\equiv 0$. \\
$(iii)$ There exists a nonnegative finite measure $d\om$ on $\dD$
such that the Herglotz representation
\begin{align}
\begin{split}
& f(z)=ic+ \oint_{\dD} d\om(\zeta) \, \f{\zeta+z}{\zeta-z}, \quad
z\in\bbD,
\\
& c=\Im(f(0)), \quad \oint_{\dD} d\om(\zeta) = \Re(f(0)) < \infty,
\end{split} \lb{C.3}
\end{align}
holds. Conversely, any function $f$ of the type \eqref{C.3} is a
Caratheodory function. \\
$(iv)$ The absolutely continuous part $d\om_{ac}$ of
$d\om$ with respect to the normalized Lebesgue measure $d\om_0$ on
$\dD$ is given by
\begin{equation}\lb{C.5}
d\om_{ac}(\ze)=\pi^{-1}\Re[f(\ze)]\,d\om_0(\ze).
\end{equation}
\end{theorem}

Next, we denote by
\begin{equation}
d\om =d\om_{\ac}+d\om_{\sc} +d\om_{\pp} \lb{C.8}
\end{equation}
the decomposition of $d\om$ into its absolutely continuous 
$({\it ac})$, singularly continuous $({\it sc})$, and pure point 
$({\it pp})$ parts with respect to Lebesgue measure on $\dD$.

\begin{theorem} $($\cite[Sects.\ 1.3, 1.4]{Si05}$)$.  \lb{tC.5}
Let $f$ be a Caratheodory function with representation \eqref{C.3}
and denote by $\La$ the set
\begin{equation}
\La=\{\ze\in\dD \st \Re[f(\ze)] \, \text{exists $($finitely or
infinitely$)$}\}. \lb{C.9}
\end{equation}
Then, $S$, $S_{\ac}$, $S_{\rm s}$, $S_{\sc}$, $S_{\pp}$ are
essential supports of $d\om$, $d\om_{\ac}$, $d\om_{\rm s}$,
$d\om_{\sc}$, $d\om_{\pp}$, respectively, where
\begin{align}
S&=\{\ze\in\La \st 0<\Re[f(\ze)]\leq\infty\}, \lb{C.10}
\\
S_{\ac}&=\{\ze\in\La \st 0<\Re[f(\ze)]<\infty\},
\lb{C.11} \\
S_{\rm s}&=\{\ze\in\La \st \Re[f(\ze)]=\infty\},
\lb{C.12} \\
S_{\sc}&=\Big\{\ze\in\La \,\Big|\, \Re[f(\ze)]=\infty, \,
\lim_{r\uparrow1}(1-r)f(r\ze)=0\Big\}, \lb{C.13} \\
S_{\pp}&=\bigg\{\ze\in\La \,\bigg|\, \Re[f(\ze)]=\infty, \,
\lim_{r\uparrow1}\left(\f{1-r}{2}\right)f(r\ze)=\om(\{\ze\})>0\bigg\}.
\lb{C.14}
\end{align}
Moreover, since
\begin{align} |\{\ze\in\dD \st |f(\ze)|=\infty\}|=0, \quad 
 \om(\{\ze\in\dD \st |f(\ze)|=\infty, \, \Re[f(\ze)]<\infty\})=0. 
\lb{C.16}
\end{align}
In addition, $S_{\rm ac}^{\prime}$, $S_{\rm s}^{\prime}$,
$S_{\sc}^{\prime}$, $S_{\pp}^{\prime}$ are essential supports of
$d\om_{\rm ac}$, $d\om_{\rm s}$, $d\om_{\sc}$, $d\om_{\pp}$,
respectively, where
\begin{align}
S_{\ac}&=\{\ze\in\La \st \Re[f(\ze)]>0\}, \lb{C.17a} \\
S_{\rm s}^{\prime}&=\{\ze\in\La \st |f(\ze)|=\infty\},
\lb{C.17} \\
S_{\sc}^{\prime}&=\Big\{\ze\in\La \,\Big|\, |f(\ze)|=\infty, \,
\lim_{r\uparrow1}(1-r)f(r\ze)=0\Big\}, \lb{C.18} \\
S_{\pp}^{\prime}&=\bigg\{\ze\in\La \,\bigg|\, |f(\ze)|=\infty, \,
\lim_{r\uparrow1}\left(\f{1-r}{2}\right)f(r\ze)=\om(\{\ze\})>0\bigg\}.
\lb{C.19}
\end{align}
\end{theorem}

Next, consider Caratheodory functions $\pm m_\pm$ of the type
\eqref{C.3},
\begin{align}
\begin{split}
& \pm m_\pm(z)=ic_\pm + \oint_{\dD} d\om_\pm(\zeta) \,
\f{\zeta+z}{\zeta-z},
\quad z\in\bbD, \lb{C.21} \\
& \,\, c_\pm\in\bbR,
\end{split}
\end{align}
and introduce the $2\times 2$ matrix-valued Caratheodory function
$M$ by
\begin{align}
&M(z)=\big(M_{j,k}(z)\big)_{j,k=0,1}, \quad z\in\bbD, \lb{C.22} \\
&M(z)=\f{1}{m_+(z)-m_-(z)}\begin{pmatrix} 1 &
\f{1}{2}[m_+(z)+m_-(z)]
\\
-\f{1}{2}[m_+(z)+m_-(z)] & -m_+(z)m_-(z) \end{pmatrix},
\lb{C.23} \\
& \hspace*{.85cm} = iC + \oint_{{\dD}} d\Om(\ze)
\frac{\ze+z}{\ze-z}, \quad z\in\bbD,  \lb{C.24} \\
& \, C=C^* = \Im[M(0)],\no
\end{align}
where $d\Om=(d\Om_{j,k})_{j,k=0,1}$ a $2\times 2$ matrix-valued
nonnegative measure satisfying
\begin{equation}
\oint_\dD d\,|\Om_{j,k}(\ze)| < \infty, \quad j,k=0,1. \lb{C.25}
\end{equation}
Moreover, we introduce the trace Caratheodory function $M^{\tr}$
\begin{align}
&M^{\tr}(z)=M_{0,0}(z)+M_{1,1}(z)=\f{1-m_+(z)m_-(z)}{m_+(z)-m_-(z)}
\lb{C.26} \\
& \hspace*{1.05cm} =ic + \oint_{\dD} d\Om^{\tr}(\ze)
\frac{\ze+z}{\ze-z}, \quad z\in\bbD, \lb{C.27}
\\ & \, c\in\bbR, \quad d\Om^{\tr}= d\Om_{0,0} +
d\Om_{1,1}. \no
\end{align}
Then,
\begin{equation}
d\Om \ll d\Om^{\tr} \ll d\Om   \lb{C.28}
\end{equation}
(where $d\mu \ll d\nu$ denotes that $d\mu$ is absolutely continuous
with respect to $d\nu$). By the Radon--Nikodym theorem, this implies that there is a 
self-adjoint integrable $2\times 2$ matrix $R(\ze)$ such
\begin{align}
d\Om(\ze) = R(\ze) d\Om^\tr(\ze).    \lb{C.29}
\end{align}
Moreover, the matrix $R(\ze)$ is nonnegative and given by
\begin{align}
R_{j,k}(\ze) = \lim_{r\uparrow1}\frac{\Re[M_{j,k}(r\ze)]}
{\Re[M_{0,0}(r\ze)+M_{1,1}(r\ze)]}  \, \text{ for a.e.\ $\ze\in\dD$}, \; j,k=0,1.
\lb{C.30}
\end{align}

Next, we identify
\begin{equation}
m_\pm(z) \, \text{ and } \, M_\pm(z,n_0), \quad z\in\bbD,
\lb{C.36}
\end{equation}
where $M_\pm(z,n_0)$ denote the half-lattice Weyl--Titchmarsh
$m$-functions defined in \eqref{3.10}--\eqref{3.11}. One then has
the following basic result (see also \cite{Si05b}).

\begin{theorem}\lb{tC.7} ${}$ \\
$(i)$ The CMV operator $U$ on $\ltz$ is unitarily equivalent to the
operator of multiplication by $I_2\id$ $($where $I_2$ is the
$2\times 2$ identity matrix and $\id(\ze)=\ze$, $\ze\in\dD$$)$ on
$L^2(\dD; d\Om(\cdot))$, and hence,
\begin{align}
\si(U) = \supp\,(d\Om) = \supp\,(d\Om^\tr),
\end{align}
where $d\Om$ and $d\Om^\tr$ are introduced in \eqref{C.24} and
\eqref{C.27}, respectively. \\
$(ii)$ The spectral multiplicity of $U$ is two if and only if
\begin{equation}
|\cM_2|>0,  \lb{C.37}
\end{equation}
where
\begin{equation}
\cM_2=\{\ze\in\La_+\,|\, m_+(\ze)\in\bbC\bs
i\bbR\}\cap\{\ze\in\La_-\,|\, m_-(\ze)\in\bbC\bs i\bbR\}.  \lb{C.38}
\end{equation}
If $|\cM_2|=0$, the spectrum of $U$ is simple. Moreover, $\cM_2$ is
a maximal set on which $U$ has uniform multiplicity two. \\
$(iii)$ A maximal set $\cM_1$ on which $U$ has uniform multiplicity
one is given by
\begin{align}
\cM_1&=\{\ze\in\La_+\cap\La_-\,|\, m_+(\ze)=
m_-(\ze)\in i\bbR\}  \no \\
&\quad \cup \{\ze\in\La_+\cap\La_-\,|\,
|m_+(\ze)|= |m_-(\ze)|=\infty\}  \no \\
&\quad \cup \{\ze\in\La_+\cap\La_-\,|\, m_+(\ze)\in i\bbR,
m_-(\ze)\in\bbC\bs i\bbR\}  \no \\
&\quad \cup \{\ze\in\La_+\cap\La_-\,|\, m_-(\ze)\in i\bbR,
m_+(\ze)\in\bbC\bs i\bbR\}. \lb{C.39}
\end{align}
In particular, $\sigma_{\rm s}(U)=\sigma_{\sc}(U)\cup
\sigma_{\pp}(U)$ is always simple.
\end{theorem}
\begin{proof}
We refer to \cite[Lemma 3.6]{GZ06} for a proof of $(i)$. To prove
$(ii)$ and $(iii)$ one observes that by $(i)$ and \eqref{C.29}, 
\begin{align}
\cN_k = \{\ze\in\si(U) \st \rank[R(\ze)]=k\}, \quad k=1,2,
\end{align}
denote the maximal sets where the spectrum of $U$ has multiplicity
one and two, respectively. Using \eqref{C.23} and \eqref{C.30} one
verifies that $\cN_k = \cM_k$, $k=1,2$.
\end{proof}

\medskip

\noindent {\bf Acknowledgments.} We are indebted to  Jonathan Breuer
for helpful discussions on this topic.


\end{document}